\documentclass[a4paper,10 pt]{amsart}
\usepackage[short]{srollens-en}[2006/10/02]
\usepackage[left = 3.5cm, right = 3.5cm, headsep = 6mm,
footskip = 10mm, top = 35mm, bottom = 35mm, footnotesep=5mm, headheight =
2cm]{geometry}

\usepackage[hypertexnames=false,
backref=page,
	pdfpagemode=UseNone,
	breaklinks=true,
	extension=pdf,
	colorlinks=true,
	linkcolor=blue,
	citecolor=blue,
	urlcolor=blue,
]{hyperref}

\title{The Fr\"olicher spectral sequence can be arbitrarily non degenerate}
\author{laura Bigalke}
\address{Laura Bigalke \\Fakult\"at f\"ur Mathematik\\Universt\"at Bielefeld\\Universit\"atsstr. 25\\33615 Bielefeld\\Germany}
\email{lbigalke@math.uni-bielefeld.de}

\author{S\"onke Rollenske}
\address{S\"onke Rollenske\\Fakult\"at f\"ur Mathematik\\Universt\"at Bielefeld\\Universit\"atsstr. 25\\33615 Bielefeld\\Germany}
\email{rollenske@math.uni-bielefeld.de}

\begin{document}

\begin{abstract}
The Fr\"olicher spectral sequence of a compact complex manifold $X$ measures the difference between Dolbeault cohomology and de Rham cohomology.

We construct for $n\geq 2$  nilmanifolds with left-invariant complex structure $X_n$ such that the $n$-th differential $d_n$ does not vanish. This replaces an earlier but incorrect example by the second author.
\end{abstract}

\maketitle

\subsection*{Introduction} 
Let $X$ be a compact complex manifold. One of the most basic invariants of the holomorphic structure of $X$ are its Hodge-numbers
\[ h^{p,q}(X)=\dim H^{p,q}(X)=\dim H^q(X, \Omega^p_X)\] 
and it is an important question how these are related to topological invariants like the Betti-numbers
\[ b_k(X)=\dim H^k(X, \IC).\]

In order to produce relations between these Fr\"olicher studied in  \cite{froelicher55} a spectral sequence connecting Dolbeault cohomology and de Rham cohomology: if we denote by $(\ka^k(X),d)$ the  complex valued de Rham complex then the decomposition of the exterior differential $d=\del +\delbar$ gives rise to a decomposition
\[\ka^k(X)=\bigoplus_{p+q=k} \ka^{p,q}(X)\]
where $\ka^{p,q}(X)$ is the space of forms of type $(p,q)$. The resulting double complex $(\ka^{p,q}(X),\del, \delbar)$ yields a spectral sequence such that
\begin{gather*}
E_0^{p,q}=\ka^{p,q}(X) \qquad d_0 = \delbar,\\
E_1^{p,q}=H^{p,q}(X) \qquad d_1 = [\del],\\
E_n^{\ast, \ast}\implies H^*_{dR}(X, \IC),
\end{gather*}
the so-called Fr\"olicher spectral sequence (see e.g. \cite{G-H}, p. 444).
It is well known that it  degenerates at the $E_1$ term for K\"ahler manifolds, which can for example be deduced from the Hodge-decomposition $H^k(X, \IC)=\bigoplus_{p+q=k} H^{p,q}(X)$.

Therefore, non-degenerating differentials in higher steps measure how much Dolbeault cohomology differs from de Rham cohomology and, in some sense, how far $X$ is from being a K\"ahler manifold.

Kodaira \cite{kodaira64} showed that for compact complex surfaces $d_1$ is always zero and the first example with $d_1\neq 0$ was the Iwasawa manifold: consider the nilpotent complex Lie group 
\[G:=\left\{ \begin{pmatrix} 1&x&y\\0&1&z\\0&0&1\end{pmatrix}\mid x,y,z\in \IC\right\}\]
and the discrete cocompat subgroup $\Gamma:=G\cap \text{Gl}(3,\IZ[i])\subset G$. Then $X:=G/\Gamma$ is a complex parallelisable nilmanifold such that $E_1\ncong E_2=E_\infty$. Later it was shown that for such manifold we have always $d_2=0$ \cite{cfg91, sakane76}.

For a long time no manifolds with $d_2\neq 0$ were known and it  was in fact speculated if $E_2=E_\infty$ holds for every compact complex manifold.

Eventually some examples with $d_2\neq 0$ were found independently by Cordero, Fern\'andez and Gray \cite{cfg87}, who used a nilmanifold with left-invariant complex structure of complex dimension 4, and Pittie \cite{pittie89}, who gave  a simply-connected example by constructing a left-invariant complex structure on $\text{Spin}(9)$.

Cordero, Fern\'andez and Gray continued their study in \cite{cfg91} finding a complex 6-dimensional nilmanifold such that $E_3\ncong E_4=E_\infty$ and together with Ugarte they showed in \cite{CFGU99} that for 3-folds several different non-degeneracy phenomena can occur up to  $E_2\ncong E_3=E_\infty$.

In \cite{rollenske07a} the second author constructed for $n\geq 2$ a series of compact complex manifolds $X_n$ and an element $[\beta_1]$ in the $E_n$-term of the Fr\"olicher spectral sequence claiming that $d_n([\beta_1])\neq 0$ (Lemma 2 in loc.cit.). This claim is incorrect: we will explain in Remark \ref{rem: false} that on the contrary $\beta_1$ induces a class in $E_\infty$.

The aim of this short note is to give a correct family of examples where with arbitrarily non-degenerate Fr\"olicher spectral sequence. This answers a  question mentioned in the book of Griffiths and Harris \cite[p.\ 444]{G-H} and repeated by Cordero, Fern\'andez and Gray.
\begin{custom}[Theorem]
For every $n\geq 2$ there exist a complex $4n-2$-dimensional compact complex manifold $X_n$ such that the Fr\"olicher spectral sequence does not degenerate at the $E_n$ term, i.e., $d_n\neq0$.
\end{custom}
We use the same method of construction as in \cite{rollenske07a} but need to introduce some extra counting variables roughly doubling the dimension of the example. We believe that in every dimension there are examples of nilmanifolds with left-invariant complex structure where the maximal possible  non-degeneracy occurs, but the structure equations might be quite complicated.
\subsection*{Construction of the example}
 \setcounter{MaxMatrixCols}{15} 
Consider the space $G_n:=\IC^{4n-2}$ with coordinates 
\[ x_1, \dots, x_{n-1}, y_1, \dots, y_n, z_1,\dots, z_{n-1}, w_1, \dots, w_n.\]
We endow $G_n$ with the structure of a real nilpotent Lie-group by identifying it with the subgroup of $\mathrm{Gl}(2n+2, \IC)$ consisting of upper triangular matrices of the form 
\begin{equation*}
A=\begin{pmatrix}
 & 1 & 0  & &&&&\dots&&0&\bar y_1&w_1\\
 &  &1 & 0 &\dots &0&-\bar z_1 & - x_1 &0  & \dots & 0 &w_2\\
&&  & \ddots& &  & & \ddots & &   &  \vdots &\vdots\\
 &&&&1 & 0 &\dots &0& -\bar z_{n-1} &- x_{n-1} & 0 & w_n\\
&&&&&1 & 0 &&\dots & & 0 & y_1\\
&&&&&&\ddots & && & \vdots & \vdots\\
\\
\\
&&&&&&&&\ddots & & \vdots & \vdots\\
&&&&&&&&& 1&0&y_n\\
&&&&&&&&& &1&z_1\\
&&&&&&&&& &&1
  \end{pmatrix}
\end{equation*}
Let $\Gamma = G_n \cap \mathrm{Gl}(2n+2, \IZ[i])$, which is a lattice in the real Lie-group $G$. 
Note that if $g\in G_n$ is a fixed element then the action on the left, $g' \mapsto gg'$ is holomorphic (with respect to the complex structure on $\IC^{4n-2}$); thus the quotient
\[ X_n = \Gamma \slash G_n\]
is a compact complex manifold; more precisely, it is a compact nilmanifold with left-invariant complex structure.

\begin{rem} \label{rem: torus bundle}
The manifold $X_n$ admits a simple geometric description in terms of principal holomorphic torus bundles:  the centre of $G_n$ is given by the matrices for which all $x_i$, $y_i$ and $z_i$ vanish and hence isomorphic (as a Lie group) to $\IC^n$. This yields an exact sequence of real Lie-groups
\[ 0\to \IC^n\to G_n\to \IC^{3n-2}\to 0\]
which is compatible with the action of $\Gamma$. Denoting by $T_k$ the quotient $\IC^k/\IZ[i]^k$ the exact sequence induces a $T_n$ principal bundle structure on  $X_n\to T_{3n-2}$.

Compare \cite{rollenske10} for a general discussion of the relation between principal holomorphic torus bundles and nilmanifolds with left-invariant complex structure and their cohomology.
\end{rem}

The space of left-invariant 1-forms $U$  is spanned by the components of $\inverse A dA$ and their complex conjugates so a basis for the forms of type $(1,0)$ is given by
\[dx_1, \dots, dx_{n-1},dy_1, \dots, dy_n, dz_1, \dots, dz_{n-1}, \omega_1,\dots, \omega_n
 \]
where
\begin{align*}
\omega_1 & = dw_1 -\bar y_1 dz_1,\\
 \omega_k &= dw_k +\bar z_{k-1}dy_{k-1} + x_{k-1} d y _k \quad(k = 2, \dots, n).
\end{align*}
For later reference we calculate the differentials of the above basis vectors:
\begin{align}
&d(dx_i)= d(dz_i) = 0 \qquad &&(i=1\dots n-1)\notag\\
&d(dy_i)= 0 \qquad &&(i=1\dots n)\notag\\
&d\omega_1 = - d\bar y_1 \wedge dz_1\notag\\
&d\omega_i = dx_{i-1}\wedge dy_i+ dy_{i-1}\wedge d\bar z_{i-1}\notag
\end{align}

The following lemma shows that the Fr\"olicher spectral sequence of $X_n$ has non-vanishing differential $d_n$ thus proving our Theorem.
\begin{lem}
The differential form $\beta_1=\bar \omega_1\wedge d\bar z_2\wedge \dots \bar dz_{n-1}$ defines a class $[\beta_1]_n\in E_n^{0,n-1}$ and
 \[ d_n([\beta_1]_n)=(-1)^{n-2}[dx_1\wedge \dots \wedge dx_{n-1}\wedge dy_n]_n\neq 0 \text{ in } E_n^{n,0}.\]
\end{lem}
\begin{proof}
By Remark \ref{rem: torus bundle} the projection to the $(x,y,z)$-coordinates endows $X_n$ with the structure of holomorphic principal torus bundle over a complex torus. By the results of \cite{con-fin01} the inclusion of left-invariant forms into the double complex $(\ka^{p,q}(X_n), \del, \delbar)$ induces an isomorphism on the $E_1$-terms of the respective spectral sequences. Thus for our purpose we may work with left-invariant forms only, that is, start with the $E_0$-term 
\[ E_0^{p,q} = \Lambda^p U\tensor \Lambda^q\bar U.\]

We say that a $(p,q)$-form $\alpha$ lives to $E_r$ if it represents a class in $E_r^{p,q}$, which is a subquotient of $E_0^{p,q}$; the resulting class will be denoted by $[\alpha]_r$.

We first show that $\beta_1$ defines a class in $E_n$. As explained in  \cite[\S14, p.161ff]{Bott-Tu}  this is equivalent to  the existence of a \emph{zig-zag} of length $n$, that is,   a collection of elements $\beta_2, \dots, \beta_{n}$ such that
\[ \beta_i\in E_0^{p+i, q-i},\quad \delbar \beta_1=0, \quad \del \beta_{i-1}+\delbar \beta_i=0 \, (i=2, \dots n).\]
We consider the following differential forms $\beta_k$ of bidegree $(k-1, n-k)$:
\begin{align*}
 \beta_2 &= \omega_2 \wedge d\bar z_2 \wedge \dots \wedge \bar z_{n-1}\\
\beta_k &= dx_1 \wedge \dots \wedge dx_{k-2} \wedge \omega_k\wedge d\bar z_k \wedge \dots \wedge \bar z_{n-1} \qquad (3\leq k\leq n-1)\\
\beta_n &= dx_1 \wedge \dots \wedge dx_{n-2} \wedge \omega_n
\end{align*}
A simple calculation shows that 
\begin{gather*}
\delbar \beta_1=0,\\
\del \beta_1=-dy_1\wedge d\bar z_1\wedge \dots\wedge d\bar z_{n-1}=-\delbar \beta_2,\\
\intertext{ and for $2\leq k\leq n-1$}
\del \beta_{k}= (-1)^{k-2} dx_1\wedge \dots \wedge dx_{k-1}\wedge dy_k\wedge d\bar z_k\wedge\dots \wedge d\bar z_{n-1}=-\delbar \beta_{k+1}.
\end{gather*}
Therefore  these elements define a zig-zag and $\beta_1$ defines a class in $E_n^{0,n-1}$.

It remains to prove that 
\[d_n[\beta_1]_n=[\del \beta_n]_n = (-1)^{n-2}[dx_1\wedge \dots \wedge dx_{n-1}\wedge dy_n]_n\]
defines a non-zero class in $E_n^{n,0}$, or equivalently, that $\beta_1$ does not live to $E_{n+1}$. In other words, we have to prove that does not exist a zig-zag of length $n+1$ for $\beta_1$. Since we are in a first quadrant double complex we have  $E_0^{n, -1}=0$ and there exists a zig-zag of length $n+1$ if and only if there exists a zig-zag $(\beta_1, \beta_2', \dots, \beta_n')$ of length $n$ such that $\del \beta_n'= 0$. 

To see that this cannot happen we put
\[U_1=\langle dx_1, \dots, dx_{n-1},dy_1, \dots, dy_n, dz_1, \dots, dz_{n-1}\rangle_\IC\text{ and }U_2 = \langle \omega_1,\dots, \omega_n\rangle_\IC,\]
 such that $U = U_1\oplus U_2$.
The above basis of $U$ and its complex conjugate induce a basis on each exterior power and a decomposition 
\[\Lambda^{n} (U\oplus \bar U) = \Lambda^n(U_1\oplus \bar U_1) \oplus S^n,\]
where $S^n$ is spanned by wedge products of basis elements, at least one of which is in $U_2\oplus \bar U_2$. 

The elements $\beta_k$ and $(-1)^{k}\del \beta_k$ are basis vectors and we decompose 
\[ E^{k-1, n-k}_0 = \beta_k\IC  \oplus V_k\text{ and } E^{k, n-k}_0 = \del\beta_k\IC  \oplus W_k,\]
where $V_k$ (resp.\ $W_k$) is spanned by all other basis elements of type $(k-1, n-k)$ (resp.\ $(k, n-k)$).
Let $\xi_k$ be the element of the dual basis such that $\xi_k\lrcorner\del\beta_k = 1$ and the contraction with any other basis element is zero.

We claim that the differentials $\del$ and $\delbar$ respect this decomposition, in the sense that
\begin{equation}\label{eq: preserve}
 \del(V_k)\subset W_k\text{ and } \delbar(V_k)\subset W_{k-1}. 
\end{equation}
More precisely, let $\alpha$ be one of the forms in our chosen basis for $\Lambda^{n-1}(U\oplus \bar U)$. Recall that $\alpha$ is a decomposable form. Then 
\begin{gather*}
\del \alpha \notin W_k\iff \xi_k \lrcorner \del \alpha \neq 0 \iff \alpha =\beta_k,\\
  \delbar \alpha \notin W_k\iff\xi_k \lrcorner \delbar \alpha \neq 0 \iff \alpha =\beta_{k+1},
\end{gather*}
which implies \eqref{eq: preserve}.
If $\alpha \in \Lambda^{n-1}(U_1\oplus \bar U_1)$ then $d \alpha = 0$ and the claim is trivial. If the form $\alpha$ contains at least two basis elements of $U_2\oplus\bar U_2$ then each summand of $d\alpha$ with respect to the basis contains at least one element of $U_2\oplus\bar U_2$, in other words, $d\alpha \in S^n$. Since $\xi_k\lrcorner S^n = 0$ the claim is true also for those elements.

The remaining elements of the basis are of the form $\pm\alpha'\wedge \omega_i$ or $\pm\alpha'\wedge \bar\omega_i$ for some $\alpha'$ in our chosen basis for $\Lambda^{n-2} (U_1\oplus \bar U_1)$. Paying special attention to the \emph{counting variable} $dy_k$ this case is easily checked by looking for solutions of the equation
\begin{multline*}
\del (\alpha'\wedge \omega_i) = (-1)^{n-1} \alpha'\wedge \del\omega_i\\
 =(-1)^{k-2} dx_1\wedge \dots \wedge dx_{k-1}\wedge dy_k\wedge d\bar z_k\wedge\dots \wedge d\bar z_{n-1}= \del\beta_k , 
\end{multline*}
and similarly in  the other cases involving either $\delbar$ or $\bar\omega_i$.

Thus if $(\beta_1, \beta_2', \dots, \beta_n')$ is any zig-zag of length $n$ for $\beta_1$ then $\beta_k' \equiv \beta_k \mod V_k$ by \eqref{eq: preserve} and, in particular,
\[ \del \beta_n' \equiv \del \beta_n \not\equiv 0 \mod W_n.\]
Thus $\beta_1$ does not live to $E_{n+1}$ and $d_n[\beta_1]_n$ is non-zero as claimed.
\end{proof}

\begin{rem}\label{rem: false}
In \cite{rollenske07a} we constructed a compact complex manifold in a very similar way and an element $[\beta_1]_n \in E_n^{0, n-1}$. However, our claim that $d_n([\beta_1]_n)\neq0$ was wrong:  while the constructed zig-zag could not be extended, the sequence of elements (in the notation of \cite[Lem.~2]{rollenske07a})
\[(\beta_1, dx_1\wedge \bar \omega_2\wedge d\bar x_3 \wedge \dots\wedge d \bar x_{n-1}, 0,0,\dots)\] gives an infinite zig-zag for $\beta_1$. In other words, the element considered  gives an element of $E_\infty$ and thus a de Rham cohomology class. 
\end{rem}

\subsection*{Acknowledgements} The second author would like to thank Daniele Grandini for pointing out the gap in \cite{rollenske07a} and the first author for showing that the gap could not be filled and for joining the search for a correct example. 
                                                                                                                                                                                                                                                                                      Both authors were supported by the DFG via the second authors Emmy-Noether project and partially via SFB 701.

%

\begin{thebibliography}{CFUG99}

\bibitem[BT82]{Bott-Tu}
Raoul Bott and Loring~W. Tu.
\newblock {\em Differential forms in algebraic topology}, volume~82 of {\em
  Graduate Texts in Mathematics}.
\newblock Springer-Verlag, New York, 1982.

\bibitem[CF01]{con-fin01}
S.~Console and A.~Fino.
\newblock Dolbeault cohomology of compact nilmanifolds.
\newblock {\em Transform. Groups}, 6(2):111--124, 2001.

\bibitem[CFG87]{cfg87}
Luis~A. Cordero, Marisa Fern{\'a}ndez, and Alfred Gray.
\newblock The {F}r{\"o}licher spectral sequence and complex compact
  nilmanifolds.
\newblock {\em C. R. Acad. Sci. Paris S{\'e}r. I Math.}, 305(17):753--756,
  1987.

\bibitem[CFG91]{cfg91}
Luis~A. Cordero, Marisa Fern{\'a}ndez, and Alfred Gray.
\newblock The {F}r{\"o}licher spectral sequence for compact nilmanifolds.
\newblock {\em Illinois J. Math.}, 35(1):56--67, 1991.

\bibitem[CFUG99]{CFGU99}
L.~A. Cordero, M.~Fernandez, L.~Ugarte, and A.~Gray.
\newblock Fr{\"o}licher spectral sequence of compact nilmanifolds with
  nilpotent complex structure.
\newblock In {\em New developments in differential geometry, Budapest 1996},
  pages 77--102. Kluwer Acad. Publ., Dordrecht, 1999.

\bibitem[Fr{\"o}55]{froelicher55}
Alfred Fr{\"o}licher.
\newblock Relations between the cohomology groups of {D}olbeault and
  topological invariants.
\newblock {\em Proc. Nat. Acad. Sci. U.S.A.}, 41:641--644, 1955.

\bibitem[GH78]{G-H}
Phillip Griffiths and Joseph Harris.
\newblock {\em Principles of algebraic geometry}.
\newblock Wiley-Interscience [John Wiley \& Sons], New York, 1978.
\newblock Pure and Applied Mathematics.

\bibitem[Kod64]{kodaira64}
K.~Kodaira.
\newblock On the structure of compact complex analytic surfaces. {I}.
\newblock {\em Amer. J. Math.}, 86:751--798, 1964.

\bibitem[Pit89]{pittie89}
Harsh~V. Pittie.
\newblock The nondegeneration of the {H}odge-de {R}ham spectral sequence.
\newblock {\em Bull. Amer. Math. Soc. (N.S.)}, 20(1):19--22, 1989.

\bibitem[Rol08]{rollenske07a}
S{\"o}nke Rollenske.
\newblock The {F}r{\"o}licher spectral sequence can be arbitrarily
  non-degenerate.
\newblock {\em Math. Ann.}, 341(3):623--628, 2008.

\bibitem[Rol11]{rollenske10}
S{\"o}nke Rollenske.
\newblock Dolbeault cohomology of nilmanifolds with left-invariant complex
  structure.
\newblock In {\em {C}omplex and {D}ifferential {G}eometry, {C}onference held at
  {L}eibniz {U}niversit\"at {H}annover, {S}eptember 14 -- 18, 2009}, volume~8
  of {\em Springer Proceedings in Mathematics}. Springer, 2011.

\bibitem[Sak76]{sakane76}
Yusuke Sakane.
\newblock On compact complex parallelisable solvmanifolds.
\newblock {\em Osaka J. Math.}, 13(1):187--212, 1976.

\end{thebibliography}
\end{document}